\documentclass[12pt]{article}
\usepackage{graphicx}
\usepackage{float,epsfig}
\usepackage[cp866]{inputenc}
\usepackage[psamsfonts]{amssymb}
\usepackage{amsfonts}
\usepackage{amsthm}

\def\Mod{\mathop{\rm Mod}\nolimits}

\begin{document}

\begin{center}\textbf{\Large S.~R.~Nasyrov}\\[2mm]
\textbf{\Large CONFORMAL MAPPINGS OF \\[1mm] STRETCHED  POLYOMINOES \\[3mm] ONTO
HALF-PLANE}
\end{center}

{\small We give an algorithm for finding conformal mappings onto the
upper half-plane and conformal modules of some types of polygons.
The polygons are obtained by stretching along the real axis
polyominoes i.~e., polygons which are connected unions of unit
squares with vertices from the integer lattice. We consider the
polyominoes of two types, so-called the $P$-pentomino and the
$L$-tetromino. The proofs are based on the Riemann-Schwarz
reflection principle and uni\-formi\-za\-tion of compact
simply-connected Riemann surfaces by rational functions.}

\section{Introduction}

Recently many authors have studied the computation of moduli of
quadrilaterals and rings and capacities of condensers and ring
domains, see, e.~g., \cite{kuhnau}, \cite{betsakos}, \cite{hakula},
\cite{dubinin}, \cite{zhang}). Special attention is devoted to the
case of polygonal boundaries of domains.

It is very important for the theory of conformal and quasiconformal
mappings to study behavior of the modules under various distortions
of domains by quasiconformal maps, for example, under the stretching
$f_H:x+iy\mapsto Hx+iy$, $H>0$, with respect to the real axis (the
problem by M.Vuorinen).

One more important direction of investigations is development of
appro\-xi\-mate methods of  calculating of conformal modules and
capacities such as method of decomposition of domain, finite element
method and others (see, e.~g., \cite{betsakos}, \cite{hakula},
\cite{zhang}, \cite{gaier_hayman1}, \cite{gaier_hayman2},
\cite{laugesen}, and \cite{falcao}).

Along with approximate methods, it is of interest obtaining exact
analytic formulas for domains of special kind. The formulas could be
used for checking results obtained by approximate methods.

One of the main model domains is so-called $L$-polygon, i.~e.,
hexagon which is the difference of two rectangles with common
vertex, in addition, one of them contains another. Such polygons
were studied by many specialists, see , e.~g., \cite{betsakos},
\cite{gaier_hayman1}, \cite{gaier_hayman2}, \cite{nas_arxiv}.  We
consider $L$-polygons of the form
\begin{equation}\label{polyomino}
P_{a_1,a}^{b_1,b}:= ([0,a]\times[0,b])\setminus
([0,a_1]\times[0,b_1]),\quad 0<a_1<a,\ 0<b_1<b.
\end{equation}
 Fixation of four
vertices turns it to a quadrilateral.  If we choose the vertices
$a_1$, $a$, $ib$, and $ib_1$, then the module of the quadrilateral
is four times less than the module of the doubly-connected domain
$([-a,a]\times[-b,b])\setminus ([-a_1,a_1]\times[-b_1,b_1])$.
Further, when we say about the module of an $L$-polygon we  will
have in mind exactly such choice of vertices.

In the paper we find conformal mappings of $L$-polygons of special
type onto the upper half-plane. The polygons are obtained from
polyominoes by stretching. We recall that polyomino is a polygon
which is the union of a finite number of squares of the same size;
the intersection of two distinct squares, if nonempty, is their
common size or vertex;  see., e.~g., \cite{golomb}. For simplicity
we consider unit squares with vertices from the integer lattice. For
$n=3$ it is said about trimino, for $n=4$ the polyominoes are called
tetrominoes, and for $n=5$ they are pentominoes.

Making use the notation (\ref{polyomino}) let us denote
$S:=P_{12}^{12}$, $P:=P_{12}^{13}$, and $R:=P_{12}^{23}$. The figure
$S$ is a trimino. In connection with their similarity of form with
Latin letters, the polygon $P$ is called a $P$-pentomino, and $R$ is
said to be an $L$-tetromino.

A conformal mapping of the trimino $S$, stretched by the map $f_H$,
onto the upper half-plane and its module are found in
\cite{borisova}, \cite{nas_arxiv}; we should note that a formula for
the module for non-stretched trimino can be obtained from results of
in~\cite{bow}, see also \cite{hakula}.

Here we find formulas for obtaining conformal mappings of the
stretched figures $P_H:=f_H(P)$  and $R_H:=f_H(R)$ onto the upper
half-plane and, as a corollary, their modules(Theorems~1 and~2).

The proofs use the Schwarz-Cristoffel symmetry principle and
uni\-for\-mi\-za\-tion of simply-connected Riemann surfaces by
rational functions.

\section{$P$-pentomino} Let $H>0$. Consider the polygon $P_H=ABCEFL$
where $A=i$, $B=3i$, $C=2H+3i$, $E=2H$, $F=H$, and $L=H+i$ (Fig.~1).
Let us denote by $\Mod(P_H)$ its module, i.~e., the extremal length
of the family of curves joining its sides $AB$ and $EF$.

Let us map the rectangle $ABCD$ onto the upper half-plane. First we
map $DMGL$ onto the part of the unit disk lying in the first quarter
of the plane $w=\varphi+i\psi$, i.~e., the domain $\{|w|<1,
\varphi>0, \psi>0\}$ by the function $f$ so that the vertices $L$,
$G$, $M$, and $D$ are mapped on $0$, $i$, $1$, and $a$; here
$a\in(0,1)$. It is easy to see that $a$ is defined from the
relations
\begin{equation}\label{lambda}
 a=\sqrt{\lambda}, \quad 2\frac{K(\lambda)}{K(\lambda')}=H.
\end{equation} Here $$K(\lambda)=\int_0^1 \frac{d\xi}{\sqrt{(1-\xi^2)(1-\lambda^2\xi^2)}}$$ is the elliptic integral of the first kind  and $\lambda'=\sqrt{1-\lambda^2}$.

By the reflection principle $f$ can be extended to the rectangle
$ABCD$. Next we extend $f$ through the segment $LD$ to the rectangle
$DEFL$. As a result, we obtain the conformal mapping of $P_H$ onto
the domain $\Sigma$ represented on Fig.~1. It is obtained by gluing
$\{|w|<1, \varphi>0, \psi<0\}$ and the upper half-plane along the
segment $[0,a]$. \vspace{1cm}

\hskip-2cm \unitlength 1mm \linethickness{0.4pt}
\ifx\plotpoint\undefined\newsavebox{\plotpoint}\fi
\begin{picture}(146.75,79.5)(0,0)
\put(32.25,19){\line(1,0){15.75}} \put(48,19){\line(0,1){52}}
\put(48,71){\line(-1,0){30.25}} \put(32.25,36){\line(0,-1){16.75}}
\put(11.75,19){\vector(1,0){53.5}} \put(17.5,14){\vector(0,1){65.5}}
\put(14.5,77){\makebox(0,0)[cc]{$y$}}
\put(64.25,15.75){\makebox(0,0)[cc]{$x$}}
\put(12.75,36.25){\makebox(0,0)[cc]{$A$}}
\put(13,71){\makebox(0,0)[cc]{$B$}}
\put(52,73.25){\makebox(0,0)[cc]{$C$}}
\put(51.75,22.75){\makebox(0,0)[cc]{$E$}}
\put(29.25,22.75){\makebox(0,0)[cc]{$F$}}
\put(34,39.25){\makebox(0,0)[cc]{$L$}}
\put(32.43,35.93){\line(1,0){.9688}}
\put(34.367,35.93){\line(1,0){.9688}}
\put(36.305,35.93){\line(1,0){.9688}}
\put(38.242,35.93){\line(1,0){.9688}}
\put(40.18,35.93){\line(1,0){.9688}}
\put(42.117,35.93){\line(1,0){.9688}}
\put(44.055,35.93){\line(1,0){.9688}}
\put(45.992,35.93){\line(1,0){.9688}}
\put(32.18,35.68){\line(0,1){.9722}}
\put(32.18,37.624){\line(0,1){.9722}}
\put(32.18,39.569){\line(0,1){.9722}}
\put(32.18,41.513){\line(0,1){.9722}}
\put(32.18,43.457){\line(0,1){.9722}}
\put(32.18,45.402){\line(0,1){.9722}}
\put(32.18,47.346){\line(0,1){.9722}}
\put(32.18,49.291){\line(0,1){.9722}}
\put(32.18,51.235){\line(0,1){.9722}}
\put(32.18,53.18){\line(0,1){.9722}}
\put(32.18,55.124){\line(0,1){.9722}}
\put(32.18,57.069){\line(0,1){.9722}}
\put(32.18,59.013){\line(0,1){.9722}}
\put(32.18,60.957){\line(0,1){.9722}}
\put(32.18,62.902){\line(0,1){.9722}}
\put(32.18,64.846){\line(0,1){.9722}}
\put(32.18,66.791){\line(0,1){.9722}}
\put(32.18,68.735){\line(0,1){.9722}}
\put(17.43,52.93){\line(1,0){.9758}}
\put(19.381,52.93){\line(1,0){.9758}}
\put(21.333,52.93){\line(1,0){.9758}}
\put(23.285,52.93){\line(1,0){.9758}}
\put(25.236,52.93){\line(1,0){.9758}}
\put(27.188,52.93){\line(1,0){.9758}}
\put(29.139,52.93){\line(1,0){.9758}}
\put(31.091,52.93){\line(1,0){.9758}}
\put(33.043,52.93){\line(1,0){.9758}}
\put(34.994,52.93){\line(1,0){.9758}}
\put(36.946,52.93){\line(1,0){.9758}}
\put(38.897,52.93){\line(1,0){.9758}}
\put(40.849,52.93){\line(1,0){.9758}}
\put(42.801,52.93){\line(1,0){.9758}}
\put(44.752,52.93){\line(1,0){.9758}}
\put(46.704,52.93){\line(1,0){.9758}}
\put(51.25,37){\makebox(0,0)[cc]{$D$}}
\put(51.5,52.75){\makebox(0,0)[cc]{$M$}}
\put(35,55.5){\makebox(0,0)[cc]{$G$}}
\put(82.5,4.75){\makebox(0,0)[cc]{Fig.~1}}
\put(75.25,41.25){\vector(1,0){73.5}} \put(32.2,19.1){\circle*{1.2}}
\put(32.2,36.){\circle*{1.2}} \put(32.2,53.){\circle*{.8}}
\put(32.2,71.){\circle*{1.2}} \put(48.0,19.1){\circle*{1.2}}
\put(48.0,36.){\circle*{1.2}} \put(48.0,53.){\circle*{1.2}}
\put(48.0,71.){\circle*{1.2}} \put(17.5,53.){\circle*{1.2}}
\put(17.5,71.){\circle*{1.2}} \put(17.5,36.){\circle*{1.2}}
\put(125.517,41.25){\line(0,1){.813}}
\put(125.497,42.063){\line(0,1){.8111}}
\put(125.438,42.874){\line(0,1){.8072}}
\multiput(125.34,43.681)(-.02748,.160306){5}{\line(0,1){.160306}}
\multiput(125.202,44.483)(-.02935,.132321){6}{\line(0,1){.132321}}
\multiput(125.026,45.277)(-.030627,.112065){7}{\line(0,1){.112065}}
\multiput(124.812,46.061)(-.031522,.096642){8}{\line(0,1){.096642}}
\multiput(124.56,46.834)(-.032152,.084444){9}{\line(0,1){.084444}}
\multiput(124.27,47.594)(-.032587,.074507){10}{\line(0,1){.074507}}
\multiput(123.944,48.339)(-.032874,.066218){11}{\line(0,1){.066218}}
\multiput(123.583,49.068)(-.033042,.059167){12}{\line(0,1){.059167}}
\multiput(123.186,49.778)(-.033113,.0530728){13}{\line(0,1){.0530728}}
\multiput(122.756,50.468)(-.0331011,.047733){14}{\line(0,1){.047733}}
\multiput(122.292,51.136)(-.0330182,.0430005){15}{\line(0,1){.0430005}}
\multiput(121.797,51.781)(-.0328728,.0387646){16}{\line(0,1){.0387646}}
\multiput(121.271,52.401)(-.0326717,.0349413){17}{\line(0,1){.0349413}}
\multiput(120.716,52.995)(-.0343275,.033316){17}{\line(-1,0){.0343275}}
\multiput(120.132,53.562)(-.0381464,.0335882){16}{\line(-1,0){.0381464}}
\multiput(119.522,54.099)(-.0397302,.0316991){16}{\line(-1,0){.0397302}}
\multiput(118.886,54.606)(-.0439685,.0317177){15}{\line(-1,0){.0439685}}
\multiput(118.227,55.082)(-.0487014,.0316591){14}{\line(-1,0){.0487014}}
\multiput(117.545,55.525)(-.0540391,.0315113){13}{\line(-1,0){.0540391}}
\multiput(116.842,55.935)(-.060129,.031259){12}{\line(-1,0){.060129}}
\multiput(116.121,56.31)(-.067171,.03088){11}{\line(-1,0){.067171}}
\multiput(115.382,56.65)(-.083832,.033717){9}{\line(-1,0){.083832}}
\multiput(114.627,56.953)(-.096039,.033314){8}{\line(-1,0){.096039}}
\multiput(113.859,57.22)(-.111476,.032706){7}{\line(-1,0){.111476}}
\multiput(113.079,57.449)(-.131752,.031806){6}{\line(-1,0){.131752}}
\multiput(112.288,57.639)(-.159767,.030457){5}{\line(-1,0){.159767}}
\put(111.489,57.792){\line(-1,0){.8053}}
\put(110.684,57.905){\line(-1,0){.8098}}
\put(109.874,57.979){\line(-1,0){.8125}}
\put(109.062,58.014){\line(-1,0){.8132}}
\put(108.249,58.01){\line(-1,0){.812}}
\put(107.437,57.965){\line(-1,0){.8089}}
\multiput(106.628,57.882)(-.20099,-.03062){4}{\line(-1,0){.20099}}
\multiput(105.824,57.76)(-.159413,-.032261){5}{\line(-1,0){.159413}}
\multiput(105.027,57.598)(-.131384,-.033293){6}{\line(-1,0){.131384}}
\multiput(104.238,57.399)(-.097212,-.029719){8}{\line(-1,0){.097212}}
\multiput(103.461,57.161)(-.085028,-.030575){9}{\line(-1,0){.085028}}
\multiput(102.695,56.886)(-.075101,-.031196){10}{\line(-1,0){.075101}}
\multiput(101.944,56.574)(-.066818,-.031637){11}{\line(-1,0){.066818}}
\multiput(101.209,56.226)(-.059772,-.031936){12}{\line(-1,0){.059772}}
\multiput(100.492,55.842)(-.0536795,-.0321201){13}{\line(-1,0){.0536795}}
\multiput(99.794,55.425)(-.0483405,-.0322075){14}{\line(-1,0){.0483405}}
\multiput(99.118,54.974)(-.0436072,-.0322126){15}{\line(-1,0){.0436072}}
\multiput(98.463,54.491)(-.0393694,-.0321461){16}{\line(-1,0){.0393694}}
\multiput(97.834,53.976)(-.035543,-.0320161){17}{\line(-1,0){.035543}}
\multiput(97.229,53.432)(-.0339488,-.0337019){17}{\line(-1,0){.0339488}}
\multiput(96.652,52.859)(-.0322748,-.0353083){17}{\line(0,-1){.0353083}}
\multiput(96.104,52.259)(-.0324326,-.0391337){16}{\line(0,-1){.0391337}}
\multiput(95.585,51.633)(-.0325301,-.0433709){15}{\line(0,-1){.0433709}}
\multiput(95.097,50.982)(-.0325596,-.0481041){14}{\line(0,-1){.0481041}}
\multiput(94.641,50.309)(-.0325111,-.0534436){13}{\line(0,-1){.0534436}}
\multiput(94.218,49.614)(-.032372,-.059537){12}{\line(0,-1){.059537}}
\multiput(93.83,48.9)(-.032124,-.066585){11}{\line(0,-1){.066585}}
\multiput(93.476,48.167)(-.031743,-.074871){10}{\line(0,-1){.074871}}
\multiput(93.159,47.419)(-.031195,-.084802){9}{\line(0,-1){.084802}}
\multiput(92.878,46.655)(-.030427,-.096992){8}{\line(0,-1){.096992}}
\multiput(92.635,45.879)(-.029359,-.112404){7}{\line(0,-1){.112404}}
\multiput(92.429,45.093)(-.033424,-.159173){5}{\line(0,-1){.159173}}
\multiput(92.262,44.297)(-.03208,-.20076){4}{\line(0,-1){.20076}}
\put(92.134,43.494){\line(0,-1){.8083}}
\put(92.045,42.685){\line(0,-1){.8117}}
\put(92.045,41.685){\line(0,-1){.8117}}
\put(108.95,17.75){\vector(0,1){59}} \put(116,40.45){\circle*{1.2}}
\put(116,45){\makebox(0,0)[cc]{$a$}}
\put(106.75,44){\makebox(0,0)[cc]{$L$}}
\put(105.25,60.75){\makebox(0,0)[cc]{$G$}}
\put(128.5,44){\makebox(0,0)[cc]{$M$}}
\put(117,37.25){\makebox(0,0)[cc]{$D$}}
\put(141.5,41.25){\circle*{1.2}}
\put(140.75,44.75){\makebox(0,0)[cc]{$1/a$}}
\put(102,41.25){\circle*{1.2}} \put(76.75,41.25){\circle*{1.2}}
\put(92.25,41.25){\circle*{1.2}} \put(109.,58.05){\circle*{1.2}}
\put(109.,41.25){\circle*{1.2}} \put(125.45,41.25){\circle*{1.2}}
\put(76,44.75){\makebox(0,0)[cc]{$-1/a$}}
\put(109.2,24.45){\circle*{1.2}} \put(125.5,39.7){\circle*{1.2}}
\thicklines \qbezier(109.3,24.5)(121.875,24.625)(125.25,39.75)
\put(115.4,41.25){\line(1,0){31.8}} \put(72,41.2){\line(1,0){37.}}
\put(116,39.75){\line(1,0){9.25}}
\put(17.75,36.25){\line(1,0){14.25}}
\put(32.25,19){\line(1,0){15.9}} \put(48,71.){\line(-1,0){30.4}}
\put(17.5,71.4){\line(0,-1){35.4}}
\put(109.,41.5){\line(0,-1){17.2}} \put(48,19.25){\line(0,1){52}}
\put(32.25,36.25){\line(0,-1){17}} \put(116,41.25){\line(0,-1){1.5}}
\put(105.25,26){\makebox(0,0)[cc]{$F$}}
\put(127.75,37.25){\makebox(0,0)[cc]{$E$}}
\put(141,37.5){\makebox(0,0)[cc]{$C$}}
\put(93,37.75){\makebox(0,0)[cc]{$N$}}
\put(14,53.5){\makebox(0,0)[cc]{$N$}}
\put(76.5,37.75){\makebox(0,0)[cc]{$B$}}
\put(31.5,74.5){\makebox(0,0)[cc]{$Q$}}
\put(116.25,72.){\makebox(0,0)[cc]{$Q(\infty)$}}
\put(100.25,45){\makebox(0,0)[cc]{$-a$}}
\put(101.5,37.75){\makebox(0,0)[cc]{$A$}}
\put(58.5,64.5){\makebox(0,0)[cc]{$z$}}
\put(58.5,64.5){\circle{3.5}} \put(127.75,67){\circle{4.2}}
\put(127.75,67){\makebox(0,0)[cc]{$w$}}
\put(147.5,39){\makebox(0,0)[cc]{$\varphi$}}
\put(105.25,75){\makebox(0,0)[cc]{$\psi$}}
\put(117.25,51){\makebox(0,0)[cc]{$\Sigma$}}
\end{picture}

Now we find a conformal mapping of $\Sigma$ onto the half-plane. To
do this we map the quarter of the unit disk $\{|\zeta|<1, \xi>0,
\eta<0\}$ in the $\zeta$-plane, $\zeta=\xi+i\eta$, onto $\Sigma$ by
the function $h$ such that the points $-i$, $0$, and $1$ are mapped
on themselves; in addition, the point $h(1)$ lies on the lower edge
of the slit along the segment $[a,1]$. Meanwhile some points $c$,
$\alpha_0\in(0,1)$ are mapped on $\infty$ and $a$; besides,
$c<\alpha_0$.

The function $h$ is extended by symmetry to a conformal map of the
unit disk onto the three-sheeted Riemann surface which is obtained
from $\Sigma$ by gluing the domains symmetric to it with respect to
the coordinate axes. It is easy to see that, in addition, points on
the unit circle are mapped on points on the unit circle, and the
extended function has a zero of the third order at the origin and
two simple poles at a pair of symmetric points $z=\pm c$.

Further we extend $h$ by symmetry to the extended complex plane. The
extended function maps the Riemann sphere onto a $5$-sheeted Riemann
surface $\mathcal{R}$. In addition, by symmetry principle, it has
zeroes at the points $\zeta=\pm 1/c$ and a pole of the third order
at the infinity. It follows that $h$ is a rational function and has
the form
$$ h(\zeta)=\zeta^3\,\frac{1-c^2\zeta^2}{\zeta^2-c^2}.$$ Because of
$\mathcal{R}$ has four branched points over the points $\pm a$, the
derivative $h'(\zeta)$ vanishes at the points $\pm \alpha_0$ and
$\pm 1/\alpha_0$.

Now we will find connections between $\alpha_0$ and $c$, and an
equation to find $\alpha_0$ if the value $a$ is fixed. We have
$h(\alpha_0)=a$, $h'(\alpha_0)=0$. Therefore, taking into account
the equality
\begin{equation}\label{log-der-h}
\frac{h'(\zeta)}{h(\zeta)}=\frac{3}{\zeta}-\frac{2\zeta}{\zeta^2-c^2}+\frac{2\zeta}{\zeta^2-1/c^2}
\end{equation}
 we obtain
\begin{equation}\label{al0}
3-\frac{2\alpha_0^2}{\alpha_0^2-c^2}+\frac{2\alpha_0^2}{\alpha_0^2-1/c^2}=0,
\end{equation}
which implies
\begin{equation}\label{alc}
\alpha_0^4-\frac{1}{3}\,(5c^2+1/c^2)\alpha_0^2+1=0.
\end{equation}
On the other hand, from (\ref{al0}) it follows that
$$
\frac{2\alpha_0^2}{\alpha_0^2-c^2}=3-\frac{2\alpha_0^2}{\alpha_0^2-1/c^2}=\frac{3-5\alpha_0^2}{1-\alpha_0^2c^2}
$$
and
$$
h(\alpha_0)=\frac{\alpha_0^3(1-\alpha_0^2c^2)}{\alpha_0^2-c^2}=\frac{\alpha_0(3-5\alpha_0^2c^2)}{2},
$$
whence, subject to $h(\alpha_0)=a$, it follows that
\begin{equation}\label{ac}
\alpha_0(3-5\alpha_0^2c^2)=2a.
\end{equation}
From (\ref{ac}) we find
\begin{equation}\label{c}
c^2=\frac{3\alpha_0-2a}{5\alpha_0^3}.
\end{equation}
Substituting (\ref{c}) to (\ref{alc}) we have
$$
\alpha_0^4-\frac{1}{3}\,\left(\frac{3\alpha_0-2a}{\alpha_0^3}+\frac{5\alpha_0^3}{3\alpha_0-2a}\right)\alpha_0^2+1=0.
$$
After simplification we obtain to the following equation to find
$\alpha_0$:
\begin{equation}\label{al0eq}
2\alpha_0^6-3a \alpha_0^5+3a\alpha_0-2a^2=0.
\end{equation}

After finding $\alpha_0$ from (\ref{al0eq}) we get the value $c$ by
(\ref{c}). Then we find the points $\alpha_1$, $\alpha_2$, and
$\alpha_3$  lying on the positive part of the abscissa axis in the
$\zeta$-plane, which correspond to the vertices $C$, $B$, and $A$ of
$P_H$. They satisfy the relations $h(\alpha_1)=1/a$,
$h(\alpha_2)=-1/a$, $h(\alpha_3)=-a$. Therefore, the points can be
got from the following equations
\begin{equation}\label{a1}
c^2x^5-x^3+(1/a)x^2-(1/a)c^2=0,
\end{equation}
\begin{equation}\label{a2}
 c^2x^5-x^3-(1/a)x^2+(1/a)c^2=0,
\end{equation}
 and
 \begin{equation}\label{a3}
 c^2x^5-x^3-ax^2+ac^2=0.
\end{equation}

After finding $\alpha_2$ and $\alpha_3$ we can get the module of
$P_H$ from the condition
\begin{equation}\label{modph}
\Mod(P_H)=\frac{K(\mu)}{K(\mu')}
\end{equation}
where
\begin{equation}\label{mu}
\mu=\frac{1+\alpha_3^2}{1-\alpha_3^2}\,\cdot\,\frac{1-\alpha_2^2}{1+\alpha_2^2};
\end{equation}
here $K(\mu)$ is the elliptic integral of the first kind, see,
e.~g., \cite{avv}.

Really, let us map $\{|\zeta|<1, \xi>0, \eta<0\}$ onto the upper
half-plane by the function
$$\omega=\left(\frac{1-w^2}{1+w^2}\right)^2.$$
The points $\alpha_k$ change under the mapping to
\begin{equation}\label{bk}
\beta_k=[(1-\alpha_k^2)/(1+\alpha_k^2)]^2,\quad 1\le k\le 3,
\end{equation} and the
points from the arc of the unit circle are mapped on points of the
negative part of the real axis. It follows that the module of $P_H$
equals the extremal length of the family of curves lying in the
upper half-plane and joining the segment $[\beta_2,\beta_3]$ with
the negative part of the real axis, i.~e., $$\frac{K(\mu)}{K(\mu')},
\quad \mu=\sqrt{\frac{\beta_2}{\beta_3}}\,.$$

It is easy to see that under the mapping the points $L$, $F$, $E$,
and $C$ change to $1$, $\infty$, $0$, and $\beta_1$.

Now we can formulate\vskip 0.3 cm

\textbf{Theorem~1.} \textit{Let $a\in(0,1)$ is defined by
$(\ref{lambda})$, $\alpha_0$ is a unique on $(0,1)$ root of the
equation $(\ref{al0eq})$, and $c$ is given by $(\ref{c})$. Let
$\alpha_1$, $\alpha_2$, and $\alpha_3$  are unique roots of the
equations $(\ref{a1})$, $(\ref{a2})$, and $(\ref{a3})$ in $(0,1)$,
and $\beta_k$, $1\le k\le 3$, are defined by $(\ref{bk})$. Then
$\Mod(P_H)$ can be calculated by $(\ref{modph})$, and the the
conformal map of the upper half-plane onto $P_H$ can be expressed
via the Schwarz-Cristoffel integral
\begin{equation}\label{sc1}
z(\zeta)=\frac{-3i}{l}\int_0^\zeta
\frac{\sqrt{t-1}}{\sqrt{t(t-\beta_1)(t-\beta_2)(t-\beta_2)}}\,dt+2H
\end{equation}
where} $$ l=\int_0^{\beta_1}
\frac{\sqrt{|t-1|}}{\sqrt{|t(t-\beta_1)(t-\beta_2)(t-\beta_2)|}}\,dt.
$$

The proof of Theorem~1 follows from the reasoning above. We only
mention some qualifying details.

The equation (\ref{al0eq}) has a unique root in $(0,1)$. Really, let
$$g(x):=2x^6-3ax^5+3ax-2a^2.$$ Then
$$g'(x)=12x^5-15ax^4+3a=12(1-a)x^5+3a(4x^5-5x^4+1)>0,\quad 0<x<1,$$ since
$a<1$ and $4x^5-5x^4+1>0$, $0<x<1$. Thus, $g$  increases on $[0,1]$.
In addition, $h(0)<0$ and $h(1)>0$, therefore, the equation $g(x)=0$
has a unique solution $\alpha_0$ in $(0,1)$.

Now we prove that $c$, defined by (\ref{c}), satisfies
$0<c<\alpha_0$. We have
$$g\left(\frac{2a}{3}\right)=-\frac{160}{243}a^6<0.$$ Since $g$
increases on $(0,1)$ and $g(\alpha_0)=0$, we have
$\frac{2a}{3}<\alpha_0$ and from (\ref{c}) we obtain $c>0$. To prove
that $c<\alpha_0$ we note that (\ref{al0eq}) can we written in the
form $\alpha_0^5(2\alpha_0-3a)+a(3\alpha_0-2a)=0$, therefore,
$a(3\alpha_0-2a)=\alpha_0^5(3a-2\alpha_0)<5a\alpha_0^4$. From the
latter inequality we obtain
$$c=\frac{3\alpha_0-2a}{5a\alpha_0^3}<\alpha_0.$$

At last, we prove that every of equations $(\ref{a1})$,
$(\ref{a2})$, and $(\ref{a3})$ has a unique solution on $(0,1)$.
From (\ref{log-der-h}), taking into account (\ref{al0}) and the
equality $h'(\alpha_0)=0$, we have
\begin{equation}\label{log-der-h1}\frac{h'(x)}{h(x)}=\frac{3(x^4-\frac{1}{3}\,
(5c^2+1/c^2)x^2+1)}{x(x^2-c^2)(x^2-1/c^2}=
\frac{3(x^2-\alpha_0^2)(x^2-1/\alpha_0^2)}{x(x^2-c^2)(x^2-1/c^2}.\end{equation}
From (\ref{log-der-h1}) we deduce that $h$ is negative and decreases
in $(0,c)$. It decreases on $(c,\alpha_0)$ from $+\infty$ to
$a\in(0,1)$ and increases on $(\alpha_0,1)$ from $a$ to $1$.
Therefore, since $0<a<1$, every of the equations $h(\alpha_1)=1/a$,
$h(\alpha_2)=-1/a$, and $h(\alpha_3)=-a$ has a unique solutions in
$(0,1)$. But they are equivalent to $(\ref{a1})$, $(\ref{a2})$, and
$(\ref{a3})$.

\vskip 0.3 cm Now we give some examples of calculating the module of
$P_H$ and the accessory parameters of the Schwarz-Cristoffel
integral (\ref{sc1}) using the package \textit{Mathematica}.

\vskip 0.5 cm

\textbf{Example~1.} For the non-stretched pentomino ($H=1$) we have
$$\lambda=3-2\sqrt{2}=0.171572875253809, \quad a=\sqrt{2}-1=0.414213562373095,$$
$$\alpha_0=           0.277046760238506,   \quad c^2=       0.025524633877222,$$
$$\alpha_1=           0.165536032447626, \quad \beta_1=     0.896160135941632,$$
$$\alpha_2=           0.154876549226231,  \quad \beta_2=    0.908495458702734,$$
$$\alpha_3=           0.138335266800084, \quad  \beta_3=    0.926301104551506,$$
$$\mu=                0.990342209151293,  \quad \Mod(P_H)=  2.137318917840447.$$

\newpage
\vskip 1 cm

\textbf{Example~2.} For $H=2$ we obtain

$$\lambda=\sqrt{2}/2=0.707106781186547, \quad      a=   0.840896415253714,$$
$$\alpha_0=          0.601898824534568, \quad     c^2=  0.113643234509673,$$
$$\alpha_1=          0.415838661746455, \quad\beta_1=   0.497227390863205,$$
$$\alpha_2=          0.301418612412185, \quad\beta_2=   0.694601063871823,$$
$$\alpha_3=          0.290931295908172, \quad\beta_3=   0.712214555130066,$$
$$\mu=               0.987557290912592, \quad \Mod(P_H)=2.056221831167256.$$

\vskip 0.2 cm

\textbf{Remark~1.} Using the method of \cite{nas_arxiv} we can
obtain the following asymp\-to\-tics: $$ \Mod(P_H)\sim{H}/{2}\,,\
H\to\infty.
$$
We see that $a$, $\lambda\to 1$, $c^2\to 1/5$, $\alpha_0$,
$\alpha_1\to 1$, $\alpha_2$, $\alpha_3\to (3-\sqrt{5})/2$,
$\beta_1\to 0$, and $\beta_2$, $\beta_3\to 1/9$ as $H\to\infty$.

For example, for $H=10$ we have $\Mod(P_H)=5.398278068084735$, and
for $H=14$ $\Mod(P_H)=7.3982087317929$.

\section{$L$-tetromino}

Consider the polygon $R_H=ABCEFL$ where $A=2i$, $B=3i$, $C=2H+3i$,
$E=2H$, $F=H$, and $L=H+2i$. Denote by $\Mod(R_H)$ its module,
i.~e., the extremal length of the family of curves joining in $R_H$
its sides $AB$ and $EF$.

Let us map conformally the rectangle $QLMC$ onto the domain
$\{|w|<1, \varphi>0, \psi>0\}$ by a function $f$ such that the
vertices $L$, $Q$, $C$, and  $M$ are mapped on $0$, $i$, $1$, and
$a$ where $a\in(0,1)$. It is easy to see that $a$ is defined from
the relations
\begin{equation}\label{a-lambda}
a=\sqrt{\lambda}, \quad 2\frac{K(\lambda)}{K(\lambda')}=H .
\end{equation}
Here, as above,  $K(\lambda)$ is the elliptic integral of the first
kind and $\lambda'=\sqrt{1-\lambda^2}$.

As in the previous section, we extend $f$ by symmetry to $R_H$, and
the extension maps conformally  $R_H$ onto the domain is obtained by
gluing the upper half of the unit disk and the quarter of the plane
$\{|w|<1, \varphi>0, \psi<0\}$ along the segment $[0,a]$. The
function $g(z)=(f(z))^2$ maps conformally $R_H$ onto the $2$-sheeted
Riemann surface $\mathcal{\widetilde{R}}$ which is a result of
gluing the unit disk, cut along $[0,1]$, and the lower half-plane
along the segment $[0,a^2]$ (Fig.~2). \vspace{1cm}

\hskip-1.6cm \unitlength 1mm \linethickness{0.4pt}
\ifx\plotpoint\undefined\newsavebox{\plotpoint}\fi
\begin{picture}(152.5,79.5)(0,0)
\put(32.25,19){\line(1,0){15.75}} \put(48,19){\line(0,1){52}}
\put(48,71){\line(-1,0){30.25}} \put(32.25,36){\line(0,-1){16.75}}
\put(11.75,19){\vector(1,0){53.5}} \put(17.5,14){\vector(0,1){65.5}}
\put(14.5,78.0){\makebox(0,0)[cc]{$y$}}
\put(64.25,15.75){\makebox(0,0)[cc]{$x$}}
\put(12.75,36.25){\makebox(0,0)[cc]{}}
\put(13,71){\makebox(0,0)[cc]{$B$}}
\put(52,73.25){\makebox(0,0)[cc]{$C$}}
\put(51.75,22.75){\makebox(0,0)[cc]{$E$}}
\put(29.25,22.75){\makebox(0,0)[cc]{$F$}}
\put(29.25,37.35){\makebox(0,0)[cc]{$G$}}
\put(34,39.25){\makebox(0,0)[cc]{}} \put(32.2,19.1){\circle*{1.2}}
\put(32.2,35.9){\circle*{1.2}} \put(32.2,53.2){\circle*{1.2}}
\put(32.2,70.8){\circle*{1.2}} \put(48.0,19.1){\circle*{1.2}}
\put(48.0,35.9){\circle*{1.2}} \put(48.0,53.2){\circle*{1.2}}
\put(48.0,70.8){\circle*{1.2}} \put(17.5,53.2){\circle*{1.2}}
\put(17.5,70.8){\circle*{1.2}} \put(32.43,35.93){\line(1,0){.9688}}
\put(34.367,35.93){\line(1,0){.9688}}
\put(36.305,35.93){\line(1,0){.9688}}
\put(38.242,35.93){\line(1,0){.9688}}
\put(40.18,35.93){\line(1,0){.9688}}
\put(42.117,35.93){\line(1,0){.9688}}
\put(44.055,35.93){\line(1,0){.9688}}
\put(45.992,35.93){\line(1,0){.9688}}
\put(32.18,35.68){\line(0,1){.9722}}
\put(32.18,37.624){\line(0,1){.9722}}
\put(32.18,39.569){\line(0,1){.9722}}
\put(32.18,41.513){\line(0,1){.9722}}
\put(32.18,43.457){\line(0,1){.9722}}
\put(32.18,45.402){\line(0,1){.9722}}
\put(32.18,47.346){\line(0,1){.9722}}
\put(32.18,49.291){\line(0,1){.9722}}
\put(32.18,51.235){\line(0,1){.9722}}
\put(32.18,53.18){\line(0,1){.9722}}
\put(32.18,55.124){\line(0,1){.9722}}
\put(32.18,57.069){\line(0,1){.9722}}
\put(32.18,59.013){\line(0,1){.9722}}
\put(32.18,60.957){\line(0,1){.9722}}
\put(32.18,62.902){\line(0,1){.9722}}
\put(32.18,64.846){\line(0,1){.9722}}
\put(32.18,66.791){\line(0,1){.9722}}
\put(32.18,68.735){\line(0,1){.9722}}
\put(51.25,37){\makebox(0,0)[cc]{$D$}}
\put(51.5,52.75){\makebox(0,0)[cc]{$M$}}
\put(35,55.5){\makebox(0,0)[cc]{$L$}}
\put(65.883,4.838){\makebox(0,0)[cc]{Fig.~2}}
\put(65.633,41.338){\vector(1,0){71.5}}
\put(90.133,44.088){\makebox(0,0)[cc]{$L$}}
\put(88.633,60.838){\makebox(0,0)[cc]{$Q$}} \thicklines
\put(88.633,26.088){\makebox(0,0)[cc]{$G$}}
\put(99.383,45.088){\makebox(0,0)[cc]{$a$}}
\put(99.383,37.838){\makebox(0,0)[cc]{$M$}}
\put(109.133,37.838){\makebox(0,0)[cc]{$D$}}
\put(111.883,44.088){\makebox(0,0)[cc]{$C$}}
\put(124.133,44.838){\makebox(0,0)[cc]{$1/a$}}
\put(124.133,37.838){\makebox(0,0)[cc]{$E$}}
\put(75.2,41.2){\circle*{1.2}} \put(85.383,41.338){\circle*{1.2}}
\put(99.5,41.9){\circle*{1.2}} \put(108.5,42.7){\circle*{1.2}}
\put(108.7,41.2){\circle*{1.2}} \put(124.883,41.2){\circle*{1.2}}
\put(99.427,42.6){\line(1,0){9.56}}
\put(101.383,18.5){\makebox(0,0)[cc]{$F(\infty)$}}
\put(76.383,37.838){\makebox(0,0)[cc]{$B$}}
\put(14,53.5){\makebox(0,0)[cc]{$A$}}
\put(76.5,37.75){\makebox(0,0)[cc]{}}
\put(31.5,74.5){\makebox(0,0)[cc]{$Q$}}
\put(83.633,45.088){\makebox(0,0)[cc]{$-a$}}
\put(84.883,37.838){\makebox(0,0)[cc]{$A$}}
\put(58.5,64.5){\makebox(0,0)[cc]{$z$}}
\put(58.5,64.5){\circle{3.5}} \put(111.133,67.088){\circle{4.2}}
\put(111.133,67.088){\makebox(0,0)[cc]{$w$}}
\put(136.5,38.7){\makebox(0,0)[cc]{$\varphi$}}
\put(88.633,75.088){\makebox(0,0)[cc]{$\psi$}}
\put(48,19.25){\line(0,1){51.5}} \put(48,19.){\line(-1,0){15.8}}
\thinlines \put(92.489,18.377){\vector(0,1){58.127}} \thicklines
\put(75.251,41.187){\line(1,0){17.238}}
\put(92.489,41.187){\line(0,-1){23.125}}
\put(99.614,41.174){\vector(1,0){38.}}
\put(99.4,41.32){\line(0,1){1.472}}
\put(92.454,41.1){\line(0,-1){23.335}}
\qbezier(92.454,58.071)(107.48,57.231)(108.718,42.365)
\qbezier(75.218,41.277)(76.5,57.143)(92.454,57.982)
\put(92.454,24.572){\circle*{1.2}} \put(92.454,41.2){\circle*{1.2}}
\put(92.454,58.1){\circle*{1.2}}
\put(17.501,53.387){\line(1,0){14.761}}
\put(17.501,70.9){\line(1,0){30.4}}
\put(32.262,53.4){\line(0,-1){34.383}}
\put(17.5,53.4){\line(0,1){17.3}} \thinlines
\put(32.103,53.228){\line(1,0){.9889}}
\put(34.08,53.217){\line(1,0){.9889}}
\put(36.058,53.206){\line(1,0){.9889}}
\put(38.036,53.195){\line(1,0){.9889}}
\put(40.014,53.184){\line(1,0){.9889}}
\put(41.991,53.173){\line(1,0){.9889}}
\put(43.969,53.162){\line(1,0){.9889}}
\put(45.947,53.151){\line(1,0){.9889}}
\end{picture}

Consider the function $h$ which maps conformally the half-disk
$\{|\zeta|<1, \eta>0\}$ onto $\mathcal{\widetilde{R}}$ such that it
keeps $0$ and maps $\pm1$ on the points with affix $1$ lying on the
different sides of the slit. Next we extend $h$ in the lower
half-disk by symmetry through the segment $[-1,1]$. As a result we
obtain the function in the unit disk which has a zero of the third
order at the point $\zeta=0$ and a simple pole at some point
$c\in(0,1)$. In addition, the points of the unit circle are mapped
to points of the unit circle. Further we extend $h$ to the extended
complex plane by symmetry. At last we conclude that $h$ is a
rational function of the form
\begin{equation}\label{rat2}
 h(\zeta)=\frac{\zeta^3(1-c\zeta)}{\zeta-c}.
\end{equation}

We may assume that $c>0$. Actually, if $c<0$, then we change
$h(\zeta)$ to
$$h(-\zeta)=\frac{\zeta^3(1+c\zeta)}{\zeta+c}$$
which has the same form as (\ref{rat2}) but correspond to the
parameter $(-c)$.

Let us find the point $\alpha_0\in (0,c)$ which is mapped on the
point $\lambda=a^2$, the endpoint of the slit of
$\mathcal{\widetilde{R}}$. We have
\begin{equation}\label{derlog}
\frac{h'(\zeta)}{h(\zeta)}=\frac{3}{\zeta}-\frac{1}{\zeta-c}-\frac{c}{1-c\zeta}.
\end{equation}
From $h'(\alpha_0)=0$ and (\ref{derlog}) it follows that
\begin{equation}\label{derloga0}
\frac{3}{\alpha_0}-\frac{1}{\alpha_0-c}-\frac{c}{1-c\alpha_0}=0
\end{equation}
or
$$
\frac{3}{\alpha_0}-\frac{1-c^2}{(\alpha_0-c)(1-c\alpha_0)}=0.$$ This
implies
\begin{equation}\label{eqa0c}
\alpha_0^2-\frac{2}{3}\,\left(2c+\frac{1}{c}\right)\alpha_0+1=0.
\end{equation}

On the other hand, from (\ref{rat2}) and (\ref{derloga0}) we have
$$
\lambda=h(\alpha_0)=\alpha_0^3\,\frac{1-
c\alpha_0}{\alpha_0-c}=3\alpha_0^2(1-c\alpha_0)-c\alpha_0^3=3\alpha_0^2-4\alpha_0^3.
$$
Therefore, $$c=\frac{3\alpha_0^2-\lambda}{4\alpha_0^3}.$$ After
substituting the obtained expression for $c$ in (\ref{eqa0c}) we get
$$
\alpha_0^2-\frac{2}{3}\,\left(\frac{3\alpha_0^2-\lambda}{2\alpha_0^3}+\frac{4\alpha_0^3}
{3\alpha_0^2-\lambda}\right)\alpha_0+1=0
$$
or, after transformation,
\begin{equation}\label{a0-2}
\alpha_0^6-3\lambda\alpha_0^4+3\lambda\alpha_0^2-\lambda^2=0.
\end{equation}

Now we find the points $\alpha_1\in (-1,0)$ and $\alpha_3\in (0,1)$
which are mapped on $\lambda$ and $1/\lambda$ under the map
(\ref{rat2}). The points can be found from the relations
$$
x^3\,\frac{1- cx}{x-c}=\lambda,
$$
$$
x^3\,\frac{1- cx}{x-c}=\frac{1}{\lambda},
$$
or, after simplification,
\begin{equation}\label{a1-2}
cx^4-x^3+\lambda x-c\lambda=0,
\end{equation}
\begin{equation}\label{a3-2}
\lambda cx^4-\lambda x^3+x-c=0.
\end{equation}

Denote $\alpha_2=c$. The desired module is equal to the extremal
length of the family of curves joining the segments $1,\alpha_1]$
and  $[\alpha_2,\alpha_3]$ in the unit upper half-disk. Let us map
the half-disk onto the upper half-plane by the function
$$\omega=\left(\frac{1+\zeta}{1-\zeta}\right)^2.$$ The points $\alpha_k$
are changed to
\begin{equation}\label{bk-2}
\beta_k=\left(\frac{1+\alpha_k}{1-\alpha_k}\right)^2.
\end{equation}
We need to calculate the extremal length of the curves joining the
segments $[-\infty,\beta_1]$ and $[\beta_2,\beta_3]$ in the upper
half-plane. It is equal to
\begin{equation}\label{mod-rh}
\Mod(R_H)=\frac{2K(\mu)}{K(\mu')}\,
\end{equation}
where
\begin{equation}\label{mu-a}
\mu=\frac{1-A}{1+A}\,,\quad
A=\sqrt{\frac{\beta_1(\beta_3-\beta_2)}{\beta_2(\beta_3-\beta_1)}}\,.
\end{equation}
\vskip 0.3cm

Thus, we have \vskip 0.3 cm

\textbf{Theorem~2.} \textit{Let $\lambda$ is defined from
$(\ref{a-lambda})$ and $\alpha_0$ is a unique solution of the
equation $(\ref{a0-2})$. Let
\begin{equation}\label{a2-2}
\alpha_2=\frac{3\alpha^2_0-\lambda}{4\alpha_0^3}\,,\end{equation}
$\alpha_1$ be a unique root of the equation $(\ref{a1-2})$ lying in
$(-1,0)$, and $\alpha_3\in (0,1)$ be a unique root of $(\ref{a3-2})$
lying in $(0,1)$. Define $\beta_k$, $1\le k\le 3$,
by~$(\ref{bk-2})$. Then $\Mod(R_H)$ can be found from
$(\ref{mod-rh})$, $(\ref{mu-a})$. The conformal map of the upper
half-plane onto $R_H$ is expressed via the Schwarz-Cristoffel
integral
\begin{equation}\label{sch-cr-2}
z(\zeta)=\frac{3i}{l}\int_{\beta_3}^\zeta
\frac{\sqrt{t-1}}{\sqrt{t(t-\beta_1)(t-\beta_2)(t-\beta_3)}}\,dt+2H\end{equation}
where }$$ l=\int_{\beta_3}^{+\infty}
\frac{\sqrt{t-1}}{\sqrt{t(t-\beta_1)(t-\beta_2)(t-\beta_3)}}\,dt.$$
 \hskip 0.3 cm

As in the case of Theorem~1, the main details of the proof are
substantiated above; we only give some additional comments for
explanation.

The equation (\ref{a0-2}) has a unique root in $(0,1)$. Really, let
$$r(x):=x^6-3\lambda x^4+3\lambda x^2-\lambda^2.$$ Then
$$r'(x)=6x(x^4-2\lambda x^2+\lambda)>0,\quad x\in (0,1),$$ since
$\lambda\in(0,1)$. Thus $h$ increases on $[0,1]$. We have
$r(0)=-\lambda^2<0$, $r(1)=1-\lambda^2>0$, therefore,  the equation
$r(x)=0$ has a unique solution $\alpha_0$ in $(0,1)$.

Now we show that $\alpha_2$, defined by (\ref{a2-2}), satisfies the
inequality $0<\alpha_2<\alpha_0$. Since $r(x)$ increases on $(0,1)$
and $r(\sqrt{\lambda/3})=-(8/27)\lambda^3<0$ we conclude that
$\alpha_0>\sqrt{\lambda/3}$. From the latter inequality it follows
that $\alpha_2>0$. From (\ref{a2-2}) we obtain
$\lambda(3\alpha_0^2-\lambda)=\alpha_0^4(3\lambda-\alpha_0^2)<4\lambda\alpha_0^4$,
therefore, $3\alpha_0^2-\lambda<4\lambda\alpha_0^4 \Rightarrow
\alpha_2<\alpha_0$.

At last, $(\ref{a1-2})$ has a unique solution  in $(-1,0)$ because
it is equivalent to the equation $r(x)=\lambda$ and $r(x)$ decreases
from $1$ to $0$ on the interval. The equation $(\ref{a3-2})$ has a
unique solution  in $(0,1)$ because it is equivalent to
$r(x)=1/\lambda$ and $r(x)$ decreases on $(0,\alpha_2)$ from
$+\infty$ to $\lambda$ and increases  on $(\alpha_2,1)$ from
$\lambda$ to $1$.\vskip 0.3 cm

Now we consider examples. \vskip 0.3 cm

\textbf{Example~3.} Non-stretched tetromino ($H=1$):

$$\lambda=3-2\sqrt{2}=   0.171572875253809,     \quad \alpha_0= 0.245789017659106,$$
$$\alpha_1=             -0.464160413208352,     \quad\beta_1=   0.133934441008549,$$
$$\alpha_2=              0.162705672169886,     \quad\beta_2=   1.928338552532678,$$
$$\alpha_3=              0.163434755042730,     \quad  \beta_3= 1.934124519781231,$$
$$A=                     0.014941124900985,     \quad \mu=      0.970557651996923,$$
$$\Mod(R_H)=             3.558625812230538.$$
\hskip 1 cm

\textbf{Example~4.} Let the coefficient of stretch  $H=2$. Then
$$\lambda=\sqrt{2}/2=0.707106781186547,        \quad \alpha_0= 0.570342096440027,$$
$$\alpha_1=         -0.872092067525734,        \quad\beta_1=   0.004668104309033,$$
$$\alpha_2=          0.362162999160609,        \quad\beta_2=   4.560775977283309,$$
$$\alpha_3=          0.401188235613446,        \quad \beta_3=  5.475355432587552,$$
$$A=                 0.013080991981584,      \quad\mu=         0.974175821903443,$$
$$\Mod(R_H)=         3.643277348370991.$$

\textbf{Remark~2.} By the methods developed in \cite{nas_arxiv} we
can obtain the following asymp\-to\-tics: $$ \Mod(R_H)\sim{H}\,,\
H\to\infty.
$$
We see that $a$, $\lambda$, $\alpha_0,\to 1$, $\alpha_1\to -1$,
$\alpha_2\to 1/2$, $\alpha_3\to 1$, $\beta_1\to 0$, and $\beta_2\to
9$, $\beta_3\to \infty$, $A\to0$ as $H\to\infty$.

For example, for $H=10$ we have $\Mod(R_H)=11.323877993993085$, and
for $H=14$ $\Mod(R_H)=15.323814523768869$.
\medskip

\textit{Acknowledgement.} The author expresses his gratitude to
Prof. M.~Vuorinen for useful comments and remarks, and to A.~Rasila
for comparing our calculations with those obtained by using the
Schwarz-Christoffel Toolbox for MATLAB created by T.~Driscoll
\cite{driscoll}, see also \cite{driscoll1}.

\vskip 0.3 cm \textsc{Institute of mathematics and mechanics, Kazan
(Volga Region) Federal University, e-mail:}
\verb"snasyrov@kpfu.ru"\vskip 0.5 cm

\end{document}